\documentclass[12pt]{article}
\usepackage{amsmath}
\usepackage{amsthm}
\usepackage{amssymb}
\usepackage{amsfonts}
\usepackage{hyperref}

\newcommand{\st}{:} 
\DeclareMathOperator{\A}{\mathit A}%
\providecommand{\majA}[1]{\sideset{^{\maj}}{_{#1}}\A}%
\providecommand{\abs}[1]{\lvert #1 \rvert}%
\DeclareMathOperator{\Des}{Des}%
\DeclareMathOperator{\des}{des}%
\DeclareMathOperator{\maj}{maj}%
\DeclareMathOperator{\nrmaj}{nrmaj}%
\DeclareMathOperator{\rmaj}{rmaj}%
\DeclareMathOperator{\rtlm}{\overleftarrow{\min}}%

\newtheorem{thm}{Theorem}[section]
\newtheorem{prop}[thm]{Proposition}
\newtheorem{lem}[thm]{Lemma}
\newtheorem{cor}[thm]{Corollary}

\theoremstyle{definition}
\newtheorem{defn}[thm]{Definition}
\newtheorem{rem}[thm]{Remark}
\newtheorem{exmp}[thm]{Example}

\title{Euler-Mahonian polynomials for $C_a\wr S_n$}

\author{Dan Bernstein\\
    Department of Mathematics\\
    The Weizmann Institute of Science\\
    Rehovot 76100, Israel\\
    \texttt{dan.bernstein@weizmann.ac.il}
}

\begin{document}

\maketitle

\section{Introduction}

Let $S_n$ denote the symmetric group on $\{1,\dots,n\}$. The
classical \emph{descent number} and \emph{major index} statistics
on $S_n$ are defined by
\[
    \des(\sigma) := \sum_{i=1}^{n-1} \chi\bigl(\sigma(i)>\sigma(i+1)\bigr)
\]
and
\[
    \maj(\sigma) := \sum_{i=1}^{n-1} i \chi\bigl(\sigma(i)>\sigma(i+1)\bigr)
\]
respectively, where $\sigma \in S_n$ and $\chi(p)$ equals $1$ if
the statement $p$ is true and $0$ if $p$ is false.

The classical Eulerian polynomials $A_n(t)$ may be defined by
\[
    A_n(t) := \sum_{\sigma\in S_n} t^{\des(\sigma)}
\]
and their ``$\maj$'' $q$-analogues, also known as the
\emph{$q$-$\maj$ Euler-Mahonian polynomials}, may be defined by
\begin{equation}\label{EQ:majAn}
    \majA{n}(t,q) := \sum_{\sigma\in S_n} t^{\des(\sigma)}
    q^{\maj(\sigma)}    .
\end{equation}
Carlitz~\cite{carlitz:qEulerian} has given a recurrence formula
for the coefficients of $\majA{n}(t,q)$: with
$\majA{n}(t,q)=\sum_{s\ge 0}t^s \majA{n,s}(q)$,
\begin{multline}\label{EQ:classicRecursion}
    \majA{n,s}(q) =\\
    (1+q+\cdots+q^s)\majA{n-1,s}(q)+(q^s+q^{s+1}+\cdots+q^{n-1})\majA{n-1,s-1}(q) .
\end{multline}
Gessel~\cite{gessel:phd} has obtained the exponential generating
function for the quotients
$\frac{\majA{n}(t,q)}{(1-t)(1-qt)\cdots(1-q^n t)}$,
\begin{equation}\label{EQ:carlitz}
    \sum_{n\ge 0}\frac{u^n}{n!}\frac{\majA{n}(t,q)}{(1-t)(1-qt)\cdots(1-q^n
    t)} = \sum_{s\ge 0}t^s e^{u(1+q+\cdots+q^s)} .
\end{equation}
Identity~\eqref{EQ:carlitz} is known as Carlitz's identity.

Let $C_a$ be the cyclic group of order $a$, and let $C_a \wr S_n$
be its wreath product with the symmetric group $S_n$, which
comprises \emph{colored permutations}.

The problem of extending the distribution of $(\des,\maj)$ to the
hyperoctahedral group $B_n = C_2 \wr S_n$ was first suggested by
Foata. Adin, Brenti and Roichman~\cite{adin:hyperoctahedral} have
given a solution to Foata's problem in the form of two pairs of
statistics, $(\mathrm{ndes},\mathrm{nmaj})$ and
$(\mathrm{fdes},\mathrm{fmaj})$. Later,
Biagioli~\cite{biagioli:Dn} has given a generalization of
Carlitz's identity to the even-signed permutation group $D_n$.
More recently, Bagno~\cite{bagno:colored} has extended the $\des$
and $\maj$ statistics to the wreath products $C_a \wr S_n$ in two
different ways, $(\mathrm{ndes},\mathrm{nmaj})$ and
$(\mathrm{ldes},\mathrm{lmaj})$, giving two further
generalizations of Carlitz's formula.

In a recent paper~\cite{regev:Wreath}, Regev and Roichman
introduced the order $<_L$ and the $C_a \wr S_n$ statistics
$\des_L$ (the $L$-descent number) and $\rtlm_L$ (number of
$L$-colored right-to-left minima), and studied the distribution of
$\rtlm_L$ on $C_a \wr S_n$ and on the subset $\{\,\sigma \in C_a
\wr S_n \st \rtlm_L(\sigma)=\des_L(\sigma) \,\}$.

Here we define the $\rmaj_{L,n}$ (\emph{$L$-reverse major index})
statistic on $C_a \wr S_n$ and study the distribution of $\des_L$
and the bi-statistic $(\des_L, \rmaj_{L,n})$. We obtain new
wreath-product analogues of the Eulerian and $q$-Euler-Mahonian
polynomials, and a generalization of Carlitz's identity (see
Corollary~\ref{CO:quotientEGF}).

\section{Preliminaries}
\subsection{The Group $C_a \wr S_n$}

Let $C_a$ be the (multiplicative) cyclic group of order $a$:
$\alpha := e^\frac{2 \pi i}{a}$ and $C_a := \{\, \alpha^t \st 0
\le t \le a-1 \,\},$ and let $C_a \wr S_n$
 be its wreath product
with $S_n$.

Elements of $C_a \wr S_n$ can be regarded as \emph{indexed
permutations} or \emph{colored permutations}---those permutations
$\sigma$ of the set $\{\, \alpha^t i \st 0\le t \le a-1, \;\; 1\le
i \le n \,\}$ satisfying $$\sigma(\beta j)=\beta\sigma(j) \quad
\forall \beta \in C_a, 1\le j \le n .$$

We shall write colored permutations using the \emph{window
notation}
\[
    \sigma = [\sigma(1),\dots,\sigma(n)]
\]
and denote
\[
    \abs{\sigma} := [\abs{\sigma(1)},\dots,\abs{\sigma(n)}] \in S_n  .
\]

Note that $S_n$ is a subgroup of $C_a \wr S_n$. $C_a \wr S_n$ is a
Coxeter group if and only if $a=1$ or $2$. In particular, for
$a=2$, $C_2 \wr S_n = B_n$ is the hyperoctahedral group, whose
elements are the \emph{signed permutations}.

\subsection{$q$-analogues}

\begin{defn}
For an integer $n\ge 1$, the $q$-analogue of $n$ is
\[
    [n]_q := 1+q+\dots+q^{n-1}.
\]
\end{defn}

\begin{defn}
For an integer $n\ge 0$, define
\[
    (\alpha;q)_n := \begin{cases}
        1,   &   \text{if $n=0$};   \\
        (1-\alpha)(1-\alpha q)\cdots(1-\alpha q^{n-1}), & \text{if
        $n>0$}.
    \end{cases}
\]
\end{defn}

\section{Statistics on $C_a \wr S_n$}

In this section we present Regev and Roichman's $<_L$ order and
various statistics based on it.

\begin{defn}[{\cite[Definition~4.4]{regev:Wreath}}]
A subset $L\subseteq \{0,1,\dots,a-1\}$ determines a linear order
$<_L$ on $\{\, \alpha^v j \st 0 \le v \le a-1,\;\; 0 \le j \le n
\,\} \cup \{0\}$ as follows:

Let $U = \{0,\dots,a-1\} \setminus L$ be the complement of $L$ in
$\{0,\dots,a-1\}$.

If $v \in L$ then $\alpha^v j <_L 0$ for every $1 \le j \le n$. If
$v\in U$ then $\alpha^v j >_L 0$ for every $1 \le j \le n$.

For $u,v \in L$ (resp. $\in U$) (not necessarily distinct) and $i
\neq j \in [n]$, $\alpha^v i <_L \alpha^u j$ if and only if $i>j$
(resp. $i<j$).

Then, for each $1 \le j \le n$, order each subset $\{\,\alpha^v j
\st v\in L\,\}$ (and each subset $\{\,\alpha^v j \st v\in U\,\}$)
in an arbitrary linear order.
\end{defn}

\begin{exmp}
Let $a=4$ and $L=\{2,3\}$, then $U=\{0,1\}$. We can choose the
following order
\begin{multline*}
\alpha^2 n <_L \alpha^3 n <_L \alpha^2(n-1) <_L \alpha^3(n-1) <_L
\dots <_L \alpha^2 <_L \alpha^3 <_L 0 \\
0 <_L \alpha <_L 1 <_L \alpha 2 <_L \dots <_L \alpha(n-1) <_L
(n-1) <_L \alpha n <_L n .
\end{multline*}
\end{exmp}

\begin{defn}[{\cite[Definition~4.6]{regev:Wreath}}]
Let $L \subseteq \{0,1,\dots,a-1\}$.

1. The \emph{$L$-descent set} of $\sigma \in C_a \wr S_n$ is
\[
    \Des_{L}(\sigma) := \{\, 0 \le i \le n-1 \st \sigma(i)>_L\sigma(i+1) \,\}
\]
where $\sigma(0):=0$.

2. The \emph{$L$-descent number} is
\[
    \des_{L}(\sigma) := \abs{\Des_{L}(\sigma)}   .
\]
\end{defn}

\begin{defn}
The \emph{$L$-reverse major index} of $\sigma \in C_a \wr S_n$ is
\[
    \rmaj_{L,n}(\sigma) := \sum_{i \in \des_L(\sigma)} n-i   .
\]
\end{defn}

\begin{rem}
In the $a=2$ case, the descent set is often defined as
\[
    \widetilde\Des_L(\sigma) := \{\, 1 \le i \le n \st \sigma(i)>_L\sigma(i+1) \,\}
\]
where $\sigma(n+1):=0$ (see for example~\cite{foata:calculI}). It
is easy to see that $\widetilde\Des_L(\sigma) = \{\,n-i \st
i\in\Des_{\{1,\dots,n\}\setminus
L}(\sigma\cdot[n,n-1,\dots,1])\,\}$, so
$\widetilde\des_L(\sigma):=\abs{\{\widetilde\Des_L(\sigma)\}} =
\des_{\{1,\dots,n\}\setminus L}(\sigma\cdot[n,\dots,1])$ and
$\widetilde\maj_L(\sigma):=\sum_{i \in \widetilde\des_L(\sigma)}i
= \rmaj_{L,n}(\sigma\cdot[n,\dots,1])$. Since multiplication by
$[n,\dots,1]$ is an involution of $C_a \wr S_n$, we get that the
bi-statistics $(\widetilde\des_L,\widetilde\maj_L)$ and
$(\des_{\{1,\dots,n\}\setminus L},\rmaj_{\{1,\dots,n\}\setminus
L,n})$ have the same distribution on $C_a \wr S_n$. Thus the
results in the following sections can be easily adapted to the
``tilde'' statistics.
\end{rem}

\begin{rem}
For $\sigma \in S_n$, let
$\rmaj_n(\sigma):=\rmaj_{\emptyset,n}(\sigma)$. A bijective
argument shows that $\maj$ and $\rmaj_n$ are equidistributed on
$\{\,\sigma \in S_n \st \des(\sigma)=s \,\}$ for every $s$. Thus
in~\eqref{EQ:majAn}, $\maj$ can be replaced by $\rmaj_n$.
\end{rem}

\begin{defn}
Define $\phi_n : C_a \wr S_{n-1} \times \{0,\dots,n-1\} \times
\{0,\dots,a-1\} \to C_a \wr S_n$ by
\[
    \phi_n(\sigma,r,t) := [\sigma_1,\dots,\sigma_r,\alpha^t n,\sigma_{r+1},\dots,\sigma_{n-1}]
\]
where $\sigma = [\sigma_1,\dots,\sigma_{n-1}]$.
\end{defn}

It is easy to see that $\phi_n$ is a bijection.

\begin{lem}
\label{LE:recursion} Let $\sigma=[\sigma_1,\dots,\sigma_{n-1}] \in
C_a \wr S_{n-1}$, $L \subseteq \{ 0,\dots,a-1 \}$ and
$\Des_L(\sigma)=\{i_1,\dots,i_s\}$, $i_1<\dots<i_s$
($s=\des_L(\sigma)$). Let $\{i_{s+1},\dots,i_n\} = \{0,\dots,n-1\}
\setminus \Des_L(\sigma)$, $n-1=i_{s+1}>\dots>i_n$ (i.e. the
non-descents of $\sigma$, from right to left). Then for $1 \le k
\le n$ and $t \in \{0,\dots,a-1\}$,
\[
    \des_L(\phi_n(\sigma,i_k,t)) = \begin{cases}
        s,      & \text{if $k<s+1$ or $k=s+1$, $t \notin L$};\\
        s+1,    & \text{if $k>s+1$ or $k=s+1$, $t \in L$} .
    \end{cases}
\]
and
\[
    \rmaj_{L,n}(\phi_n(\sigma,i_k,t)) = \begin{cases}
        \rmaj_{L,n-1}(\sigma)+k,   & \text{if $t \in L$};\\
        \rmaj_{L,n-1}(\sigma)+k-1, & \text{if $t \notin L$},
    \end{cases}
\]
\end{lem}

\begin{proof}
We consider the three possible cases:

\noindent{\bf Case 1.} $1\le k\le s$. In this case,
\[
    \tilde\sigma:=\phi_n(\sigma,i_k,t) = [\sigma_1,\dots,\sigma_{i_1},\dots,\sigma_{i_k},\alpha^t
    n,\sigma_{i_k+1},\dots,\sigma_{i_{k+1}},\dots,\sigma_{i_s},\dots,\sigma_{n-1}],
\]
thus the descents to the right of $\sigma_{i_k+1}$ are shifted one
position to the right, and the $k-1$ descents to the left of
$\sigma_{i_k}$ remain in place. If $t \in L$, then
$\sigma_{i_k}>_L \alpha^t n <_L \sigma_{i_k+1}$, so the descent at
$i_k$ is also shifted one position to the right. If $t \notin L$,
then $\sigma_{i_k}<_L \alpha^t n >_L \sigma_{i_k+1}$, leaving the
descent at $i_k$. The contribution to $\rmaj_{L,n}(\tilde\sigma)$
of each descent shifted one position to the right is the same as
its contribution to $\rmaj_{L,n-1}(\sigma)$, whereas the
contribution to $\rmaj_{L,n}(\tilde\sigma)$ of each descent left
in place is $1$ more than its contribution to
$\rmaj_{L,n-1}(\sigma)$.

\noindent{\bf Case 2.} $k=s+1$. In this case,
\[
    \tilde\sigma:=\phi_n(\sigma,n-1,t) = [\sigma_1,\dots,\sigma_{n-1},\alpha^t n],
\]
thus all $s=k-1$ descents remain in place, each contributing $1$
more to $\rmaj_{L,n}(\tilde\sigma)$ than to
$\rmaj_{L,n-1}(\sigma)$. If and only if $t\in L$, $\sigma_{n-1}>_L
\alpha^t n$ so there is an additional descent at $n-1$, which
contributes $1$ to $\rmaj_{L,n}(\tilde\sigma)$.

\noindent{\bf Case 3.} $s+1<k\le n$. In this case,
\[
    \tilde\sigma:=\phi_n(\sigma,i_k,t) = [\sigma_1,\dots,\sigma_{i_r},\dots,\sigma_{i_k},\alpha^t n,
    \sigma_{i_k+1},\dots,\sigma_{i_{r+1}},\dots,\sigma_{n-1}]
\]
where $r$ is the number of descents to the left of $\sigma_{i_k}$,
whence \[i_k = (n-1)-(k-(s+1))-(s-r) = n-k+r .\] In this case the
descents to the right of $\sigma_{i_k}$ are shifted one position
to the right, and the $r$ descents to the left of $\sigma_{i_k}$
remain in place. The contribution to $\rmaj_{L,n}(\tilde\sigma)$
of each descent shifted one position to the right is the same as
its contribution to $\rmaj_{L,n-1}(\sigma)$, whereas the
contribution to $\rmaj_{L,n}(\tilde\sigma)$ of each descent left
in place is $1$ more than its contribution to
$\rmaj_{L,n-1}(\sigma)$. If $t\in L$, then there is an additional
descent at $i_k=n-k+r$, whose contribution to
$\nrmaj_{L,n}(\tilde\sigma)$ is $k-r$. If $t \notin L$, then the
additional descent is at $i_k+1 = n-k+r+1$, contributing $k-r-1$
to $\nrmaj_{L.n}(\tilde\sigma)$.
\end{proof}

\section{$C_a \wr S_n$ $q$-$\maj$ Euler-Mahonian
Polynomials}\label{SEC:EulerMahonian}

In this section we define $q$-$\maj$ Euler-Mahonian polynomials
for $C_a \wr S_n$ and give generalizations of the results by
Carlitz and Gessel.

For $L\subseteq \{ 0,\dots,a-1 \}$, let $\majA{a,L,n}(t,q)$ be the
generating polynomial for $C_a \wr S_n$ by the bi-statistic
$(\des_L,\rmaj_{L,n})$, i.e.
\[
    \majA{a,L,n}(t,q) := \sum_{\sigma \in C_a \wr S_n}
    t^{\des_L(\sigma)} q^{\rmaj_{L,n}(\sigma)}  .
\]

\begin{rem}\label{RE:indep}
For $n=1$,
\[
\majA{a,L,1}(t,q) = \sum_{t=0}^{a-1}
t^{\des_L([\alpha^t])}q^{\rmaj_{L,1}([\alpha^t])}  = \ell t q
+(a-\ell)
\]
where $\ell=\abs{L}$ depends only on the number of elements in $L$
and not on the choice of elements.
\end{rem}

The following is a generalization of~\eqref{EQ:classicRecursion}.
\begin{prop}\label{PR:recursion}
With $\majA{a,L,n}(t,q) = \sum_{s \ge 0} t^s \majA{a,L,n,s}(q)$,
the coefficients $\majA{a,L,n,s}(q)$ satisfy the recurrence
\begin{multline}\label{EQ:recursion}
    \majA{a,L,n,s}(q) =\\ (a[s+1]_q-\ell)\majA{a,L,n-1,s}(q)+(a q^s [n-s]_q+\ell
    q^n)\majA{a,L,n-1,s-1}(q)
\end{multline}
where $\ell=\abs{L}$.
\end{prop}
\begin{proof}
By definition,
\[
\majA{a,L,n,s}(q) = \sum_{\substack{\tilde\sigma \in C_a \wr
S_n\\\des_L(\tilde\sigma)=s}}q^{\rmaj_{L,n}(\tilde\sigma)}   .
\]
By the bijectivity of $\phi_n$ and Lemma~\ref{LE:recursion},
\[
    \{\,\tilde\sigma \in C_a \wr S_n \st \des_L(\tilde\sigma) = s\,\}
    = \phi_n(A \uplus B \uplus C \uplus D)
\]
where $\uplus$ denotes disjoint union and
\begin{multline*}
    A := \{\,(\sigma,i_k,t) \st \sigma \in C_a \wr S_{n-1},\;\;
     \Des_L(\sigma)=\{i_1,\dots,i_s\},\;\; 1\le k\le s,\;\; t\in L
    \,\}
\end{multline*}
\begin{multline*}
    B := \{\,(\sigma,i_k,t) \st \sigma \in C_a \wr S_{n-1},\;\; \\
     \{0,\dots,n-1\}\setminus \Des_L(\sigma)=\{i_s,\dots,i_n\},\;\; s\le k\le n,\;\; t\in L \,\}
\end{multline*}
\begin{multline*}
    C := \{\,(\sigma,i_k,t) \st \sigma \in C_a \wr S_{n-1},\;\;\\
    \Des_L(\sigma)\uplus\{n-1\}=\{i_1,\;\;\dots,i_{s+1}\},\;\; 1\le k\le s+1,\;\; t\notin
    L\, \}
\end{multline*}
\begin{multline*}
    D := \{\,(\sigma,i_k,t) \st \sigma \in C_a \wr S_{n-1},\;\;\\
        \{0,\dots,n-2\}\setminus\Des_L(\sigma)=\{i_{s+1},\;\;\dots,i_n\},\;\; s+1\le k\le n,\;\; t\notin L
    \,\}  .
\end{multline*}
Note that in the definition of $A$ and $C$, $\des_L(\sigma)=s$,
whereas in the definition of $B$ and $D$, $\des_L(\sigma)=s-1$.

By the second part of Lemma~\ref{LE:recursion},
\[
\begin{split}
\sum_{\substack{\sigma \in C_a \wr
S_n\\\des_L(\sigma)=s}}q^{\rmaj_{L,n}(\sigma)}
    &=  \sum_{(\sigma,i_k,t) \in A}q^{\rmaj_{L,n}(\phi_n(\sigma,i_k,t))}
    +  \sum_{(\sigma,i_k,t) \in B}q^{\rmaj_{L,n}(\phi_n(\sigma,i_k,t))} \\
    &\quad + \sum_{(\sigma,i_k,t) \in C}q^{\rmaj_{L,n}(\phi_n(\sigma,i_k,t))}
    +  \sum_{(\sigma,i_k,t) \in D}q^{\rmaj_{L,n}(\phi_n(\sigma,i_k,t))} \\
    &= \sum_{(\sigma,i_k,t) \in A}q^{\rmaj_{L,n-1}(\sigma)+k}
    +  \sum_{(\sigma,i_k,t) \in B}q^{\rmaj_{L,n-1}(\sigma)+k} \\
    &\quad + \sum_{(\sigma,i_k,t) \in C}q^{\rmaj_{L,n-1}(\sigma)+k-1}
    +  \sum_{(\sigma,i_k,t) \in D}q^{\rmaj_{L,n-1}(\sigma)+k-1} \\
    &= \ell q[s]_q\sum_{\substack{\sigma \in C_a\wr S_{n-1}\\\des_L(\sigma)=s}}q^{\rmaj_{L,n-1}(\sigma)}\\
    &\quad + \ell q^s[n-s+1]_q\sum_{\substack{\sigma \in C_a\wr S_{n-1}\\\des_L(\sigma)=s-1}}
            q^{\rmaj_{L,n-1}(\sigma)}\\
&\quad + (a-\ell) [s+1]_q\sum_{\substack{\sigma \in C_a\wr
S_{n-1}\\\des_L(\sigma)=s}}
            q^{\rmaj_{L,n-1}(\sigma)}\\
&\quad + (a-\ell) q^s[n-s]_q\sum_{\substack{\sigma \in C_a\wr
S_{n-1}\\\des_L(\sigma)=s-1}}q^{\rmaj_{L,n-1}(\sigma)}\\
&= (a[s+1]_q-\ell)\majA{a,L,n-1,s}\\
&\qquad + (a q^s[n-s]_q+\ell q^n)\majA{a,L,n-1,s-1}    .\qedhere
\end{split}
\]
\end{proof}

By Proposition~\ref{PR:recursion} and Remark~\ref{RE:indep},
$\majA{a,L,n}(t,q)$ does not depend on the choice of elements in
$L$ but only on their number $\ell=\abs{L}$. Therefore the
polynomials
\begin{align*}
    \majA{a,\ell,n}(t,q)&:=\majA{a,L,n}(t,q) \qquad \abs{L}=\ell\\
\intertext{and}
    \majA{a,\ell,n,s}(q)&:=\majA{a,L,n,s}(q) \qquad \abs{L}=\ell
\end{align*}
are well-defined.

The following lemma is a generalization
of~\cite[equation~(10.3)]{foata:qSeriesRev2}.
\begin{lem}
For every three integers $n>0$, $a>0$, $0\le \ell \le a$,
\begin{multline}\label{EQ:recursion2}
    (1-q)\majA{a,\ell,n}(t,q) =\\
    (a-(1-q)\ell)(1-tq^n)\majA{a,\ell,n-1}(t,q)-aq(1-t)\majA{a,\ell,n-1}(tq,q)
    .
\end{multline}
\end{lem}
\begin{proof}
Multiply both sides of~\eqref{EQ:recursion} by $(1-q)t^s$ and sum
over all $-\infty<s<\infty$ to get
\[
\begin{split}
&(1-q)\majA{a,\ell,n}(t,q)\\
    &= \sum_{s=-\infty}^\infty(a(1-q^{s+1})-(1-q)\ell)t^s\majA{a,\ell,n-1,s}(q) \\
    &\quad +\sum_{s=-\infty}^\infty(a q^s(1-q^{n-s})+(1-q)\ell q^n)t^s\majA{a,\ell,n-1,s-1}(q)\\
    &= \sum_{s=-\infty}^\infty\bigl((a-(1-q)\ell)t^s-aq(qt)^s\bigr)\majA{a,\ell,n-1,s}(q) \\
    &\quad +\sum_{s=-\infty}^\infty t\bigl((-a q^n+(1-q)\ell q^n)t^s+a(qt)^s\bigr)\majA{a,\ell,n-1,s}(q)\\
    &= \sum_{s=-\infty}^\infty\bigl((a-(1-q)\ell)(1-tq^n)t^s-aq(1-t)(qt)^s\bigr)\majA{a,\ell,n-1,s}(q) \\
    &=
    (a-(1-q)\ell)(1-tq^n)\majA{a,\ell,n-1}(t,q)-aq(1-t)\majA{a,\ell,n-1}(tq,q)
    .\qedhere
\end{split}
\]
\end{proof}

\begin{prop}\label{PR:quotient}
For every three integers $n>0$, $a>0$, $0\le \ell \le a$,
\begin{equation}\label{EQ:quotient}
    \frac{\majA{a,\ell,n}(t,q)}{(t;q)_{n+1}} = \sum_{s\ge 0} t^s
    (a[s+1]_q-\ell)^n   .
\end{equation}
\end{prop}
\begin{proof}
By induction on $n$. For $n=1$,
\[
    \frac{\majA{a,\ell,1}(t,q)}{(t;q)_{2}} = \frac{a -\ell (1-tq)}{(1-t)(1-tq)} =
    \sum_{s\ge 0}t^s\left( a[s+1]_q -\ell \right)   .
\]
For $n>1$, divide both sides of~\eqref{EQ:recursion2} by
$(1-q)(t;q)_{n+1}$ and use the induction hypothesis to get
\[
\begin{split}
&\frac{\majA{a,\ell,n}(t,q)}{(t;q)_{n+1}}\\
    &= \frac{a-(1-q)\ell}{(1-q)(t;q)_n}\majA{a,\ell,n-1}(t,q) -
    \frac{aq(1-t)}{(1-q)(t;q)_{n+1}}\majA{a,\ell,n-1}(tq,q)\\
    &= \left(\frac{a}{1-q}-\ell\right)\sum_{s\ge
    0}t^s(a[s+1]_q-\ell)^{n-1} - \frac{aq}{1-q}\sum_{s\ge
    0}t^s q^s (a[s+1]_q-\ell)^{n-1}\\
    &= \sum_{s\ge 0}t^s\left(\frac{a}{1-q}-\ell-\frac{a
    q^{s+1}}{1-q}\right)(a[s+1]_q-\ell)^{n-1}\\
    &=\sum_{s\ge 0}t^s(a[s+1]_q-\ell)^n . \qedhere
\end{split}
\]
\end{proof}

As a corollary, we get a generalization of~\eqref{EQ:carlitz}.
\begin{cor}\label{CO:quotientEGF}
For every $a>0$, $0\le \ell \le a$,
\begin{equation}\label{EQ:quotientEGF}
    \sum_{n\ge 0}\frac{u^n}{n!}\frac{\majA{a,\ell,n}(t,q)}{(t;q)_{n+1}} = \sum_{s\ge 0} t^s
    e^{u(a[s+1]_q-\ell)}   .
\end{equation}
\end{cor}

\section{$C_a\wr S_n$ Eulerian Polynomials}

For $L\subseteq \{ 0,\dots,a-1 \}$, let $A_{a,L,n}(t)$ be the
generating polynomial for $C_a \wr S_n$ by the statistic $\des_L$,
i.e.
\[
    A_{a,L,n}(t) := \sum_{\sigma \in C_a \wr S_n}
    t^{\des_L(\sigma)}   .
\]
Clearly $A_{a,L,n}(t) = \majA{a,L,n}(t,1)$ and the polynomials
\[
    A_{a,\ell,n}(t) := A_{a,L,n}(t) \qquad \abs{L}=\ell
\]
are well-defined.

\begin{prop}
\begin{enumerate}
\item
With $A_{a,\ell,n}(t)=\sum_{s\ge 0}t^s A_{a,\ell,n,s}$, the
coefficients $A_{a,\ell,n,s}$ satisfy the recurrence
\begin{equation}\label{EQ:desRecursion}
    A_{a,\ell,n,s} =
    (a(s+1)-\ell)A_{a,\ell,n-1,s}+(a(n-s)+\ell)A_{a,\ell,n-1,s-1}
    .
\end{equation}
\item
\begin{equation}\label{EQ:desRecursion2}
    A_{n,a,\ell}(t) =
    \bigl(a-\ell+(a(n-1)+\ell)t\bigr)A_{a,\ell,n-1}(t)+at(1-t)A'_{a,\ell,n-1}(t).
\end{equation}
\item
\begin{equation}\label{EQ:desQuotient}
    \frac{A_{a,\ell,n}(t)}{(1-t)^{n+1}} = \sum_{s\ge
    0}t^s(a(s+1)-\ell)^n.
\end{equation}
\item
\begin{equation}\label{EQ:desQuotientEGF}
    \sum_{n\ge 0}\frac{u^n}{n!}\frac{A_{a,\ell,n}(t)}{(1-t)^{n+1}} = \sum_{s\ge
    0}t^s e^{u(a(s+1)-\ell)} = \frac{e^{(a-\ell)u}}{1-te^{au}}.
\end{equation}
\item
\begin{equation}\label{EQ:desQuotientGF}
    \sum_{n\ge 0}\frac{u^n}{n!} A_{a,\ell,n}(t) = \frac{1-t}{-te^{\ell(1-t)u}+e^{(a-\ell)u(t-1)}}.
\end{equation}
\end{enumerate}
\end{prop}

\begin{proof}
\eqref{EQ:desRecursion},~\eqref{EQ:desQuotient}
and~\eqref{EQ:desQuotientEGF} follow
from~\eqref{EQ:recursion},~\eqref{EQ:quotient}
and~\eqref{EQ:quotientEGF} respectively by setting $q=1$.
\eqref{EQ:desQuotientGF} follows from~\eqref{EQ:desQuotientEGF} by
substituting $(1-t)u$ for $u$ and multiplying both sides by $1=t$.
Finally, to see~\eqref{EQ:desRecursion2}, multiply both sides
of~\eqref{EQ:desRecursion} by $t^s$ and take the sum over all
$s\ge 0$ to get
\[
\begin{split}
A_{a,\ell,n}(t)
    &= \sum_{s\ge 0}(a(s+1)-\ell)t^s A_{a,\ell,n-1,s} +
    \sum_{s\ge0}(a(n-s)+\ell)t^s A_{a,\ell,n-1,s-1} \\
    &= \sum_{s\ge 0}(a(s+1)-\ell)t^s A_{a,\ell,n-1,s} +
    t \sum_{s\ge0}(a(n-s-1)+\ell)t^s A_{a,\ell,n-1,s} \\
    &= \sum_{s\ge 0}\bigl(a-\ell+(a(n-1)+\ell)t\bigr)t^s A_{a,\ell,n-1,s} \\
    &\quad + t \sum_{s\ge 0}a(1-t) s t^{s-1} A_{a,\ell,n-1,s} \\
    &=\bigl(a-\ell+(a(n-1)+\ell)t\bigr)A_{a,\ell,n-1}(t)+at(1-t)A'_{a,\ell,n-1}(t) . \qedhere
\end{split}
\]
\end{proof}

\end{document}